\documentclass[11pt,a4paper]{article}
\usepackage[utf8x]{inputenc}
\usepackage[english]{babel}
\usepackage[T1]{fontenc}
\usepackage{amsmath}
\usepackage{amsfonts}
\usepackage{amssymb}
\usepackage{makeidx}
\usepackage{graphicx}
\usepackage[left=3cm,right=2cm,top=3cm,bottom=2cm]{geometry}
\usepackage{setspace}
\usepackage{amsthm}
\usepackage{mathtools}
\usepackage{cite}

\usepackage[colorlinks=true,urlcolor=blue,
citecolor=red,linkcolor=blue,linktocpage,pdfpagelabels,
bookmarksnumbered,bookmarksopen]{hyperref}

\addto\captionsenglish{}

\newtheorem{definition}{Definition}[section]

\newtheorem{proposition}{Proposition}[section]

\newtheorem{theorem}{Theorem}[section]

\newtheorem{example}{Example}[section]

\newtheorem{lemma}{Lemma}[section]

\newtheorem{observation}{Remark}[section]

\newtheorem{corollary}{Corolary}[section]

\numberwithin{equation}{section}

\newcommand{\bo}{\begin{observation}}
\newcommand{\eo}{\end{observation}}
\newcommand{\bd}{\begin{definition}}
\newcommand{\ed}{\end{definition}}
\newcommand{\bp}{\begin{proposition}}
\newcommand{\ep}{\end{proposition}}
\newcommand{\bt}{\begin{theorem}}
\newcommand{\et}{\end{theorem}}
\newcommand{\bc}{\begin{corollary}}
\newcommand{\ec}{\end{corollary}}
\newcommand{\bl}{\begin{lemma}}
\newcommand{\el}{\end{lemma}}
\newcommand{\be}{\begin{example}}
\newcommand{\ee}{\end{example}}
\newcommand{\beq}{\begin{equation}}
\newcommand{\eeq}{\end{equation}}
\newcommand{\beqa}{\begin{equation*}}
\newcommand{\eeqa}{\end{equation*}}
\newcommand{\R}{\mathbb{R}}
\newcommand{\RN}{\mathbb{R}^{N}}
\newcommand{\N}{\mathbb{N}}
\newcommand{\WumN}{ W^{1,N}(\mathbb{R}^{N}) }
\newcommand{\Wump}{ W^{1,p}(\mathbb{R}^{N}) }

\newcommand{\LN}{ L^{N}(\mathbb{R}^N) }
\newcommand{\Lp}{ L^{p}(\mathbb{R}^N) }

\newcommand{\Lt}{L^{t}(\mathbb{R}^N)}
\newcommand{\B}{\mathcal{B}}
\newcommand{\A}{\mathcal{A}}
\newcommand{\F}{\mathcal{F}}
\newcommand{\K}{\mathcal{K}}

\newcommand{\intR}{\displaystyle\int\limits_{\mathbb{R}^N}}
\newcommand{\un}{u_{n}}

\newcommand{\vn}{v_{n}}

\newcommand{\yn}{y_{n}}

\newcommand{\until}{\tilde{u}_n}

\newcommand{\RA}{\rightarrow}
\newcommand{\IC}{\hookrightarrow}
\newcommand{\CF}{\rightharpoonup}
\newcommand{\ds}{\displaystyle\int\limits}
\newcommand{\intRR}{\iint\limits_{\RN \times \RN}}

\usepackage{times, color, xcolor}
\addto\captionsportuguese{
}

\begin{document}

 	\title{$ (p, N)-$Choquard logarithmic equation involving a nonlinearity with exponential critical growth: existence and multiplicity
		\thanks{The first author was supported    by  Coordination of Superior Level Staff Improvement-(CAPES) -Finance Code 001 and  S\~ao Paulo Research Foundation- (FAPESP), grant $\sharp $ 2019/22531-4,
while the second  author was supported by  National Council for Scientific and Technological Development -(CNPq),   grant $\sharp $ 307061/2018-3 and FAPESP  grant $\sharp $ 2019/24901-3.
			}}
	\author{
		Eduardo  de S. Böer \thanks{Corresponding author} \thanks{ E-mail address: eduardoboer04@gmail.com Tel. +55.51.993673377}  and Ol\'{\i}mpio H. Miyagaki  \footnote{ E-mail address: ohmiyagaki@gmail.com, Tel.: +55.16.33519178 (UFSCar).}\\
		{\footnotesize Department of Mathematics, Federal University of S\~ao Carlos,}\\
		{\footnotesize 13565-905 S\~ao Carlos, SP - Brazil}\\ }
\noindent
				
	\maketitle

\noindent \textbf{Abstract:} The present work is concerned with the following version of Choquard Logarithmic equations $ -\Delta_p u -\Delta_N u + a|u|^{p-2}u + b|u|^{N-2}u + \lambda (\ln|\cdot|\ast G(u))g(u) = f(u) \textrm{ \ in \ } \mathbb{R}^N $ , where $ a, b, \lambda >0 $, $ \max\{\frac{N}{2}, 2 \} < p< N $, $f, g: \R \RA \R $ are continuous functions that behave like $ \exp(\alpha |u|^{\frac{N}{N-1}}) $ at infinity,  for $ \alpha >0 $, and that has polynomial growth, respectively, and $ G(s)=\int\limits_{0}^{s}g(\tau)d\tau $. We prove the existence of a nontrivial solution at the mountain pass level and a nontrivial ground state solution. Also, using a version of the Symmetric Mountain-Pass Theorem, we get infinitely many solutions.

\vspace{0.5 cm}


\section{Introduction}

In the present paper we are interested in the existence and multiplicity of solutions for $ (p, N)-$Choquard logarithmic equations, with a nonlinearity of Moser-Trudinger type, 
\begin{equation} \label{P}
-\Delta_p u -\Delta_N u + a|u|^{p-2}u + b|u|^{N-2}u + \lambda (\ln|\cdot|\ast G(u))g(u) = f(u) \textrm{ \ in \ } \mathbb{R}^N,
\end{equation}
with $ a \equiv 1 $, $ \lambda \equiv 1 $, $ N \geq 3 $, $ \max\{2, \frac{N}{2}\} < p < N $, $ g:\R \RA \R $  a continuous function with polynomial growth and $f: \R \RA \R $ a continuous function that behaves like $ \exp(\alpha |u|^{N'}) $ at infinity, where $ N' =\frac{N}{N-1} $ and whose primitives are given by $ F(s)=\int\limits_{0}^{s}f(\tau)d\tau $ and $ G(s)=\int\limits_{0}^{s}g(\tau)d\tau $, respectively.

We begin providing a quick overview of literature, focusing especially in works involving the logarithmic kernel.  The reader also should pay attention that some topics have a really vast quantity of results, so we choose only some of them to mention in our introduction. We start mentioning that Choquard equations have been extensively studied, without $N$-Laplacian operator, in the case $ N\geq 3 $, $ p=2 $ and $ f(u)=b|u|^{p-2}u $, it occurred mostly because such equations have many important applications in physics. Also, the kernel in this case presents an easier behaviour to deal with. We can cite as examples \cite{[12], [22], [20] , [16], [18]}, in which the authors treated many variations of Choquard equation.

On the other hand, when dealing with logarithmic forms of Choquard equations, there are few references addressing the topic. We can cite some recent works \cite{[6], [cjj], [10], [4]} where the authors have established the results on existence and multiplicity of solutions for equations with polynomial nonlinearity. Moreover, in \cite{[alves]}, the authors have dealt with equation (\ref{P}) considering a nonlinearity with exponential critical growth. They have proved the existence of a ground state. Finally, in \cite{[boer3]}, the authors have proved existence and multiplicity results for the $ p-$fractional Laplacian operator.  We would like to say that, while dealing with exponential nonlinearities, some additional difficulties arise what makes necessary to consider some stronger conditions over the nonlinearity. 

To finish this overview, we cite some authors that study PDEs involving $ N-$Laplacian with exponential critical nonlinearities, e.g. \cite{[lam], [doO], [zhang]}, and some that dealt with $ (p, q) $-Laplacian equations with  exponential behaviour problems, for instance, \cite{[pq1], [pq2], [pq3], [pq4]}. It is important to say that, in order to work with an exponential nonlinearity, we need to take $ q =N $ in the operator, where $ N $ stands for the dimension of the euclidian space.  The work that inspired us to study such problem was \cite{[pucci]}.

The present work intend to extend or complement the results found in the above mentioned papers, combining the critical  exponential growth with the $(p, N)-$Laplacian logarithmic Choquard equation. Moreover, it is the first result considering a more general nonlinearity in the convolution with the logarithm, which increases significantly the difficulties to obtain the results and the scope of equations solved.

A great difficulty in our work, as mentioned in the previous works that deal with Choquard logarithmic equations is to guarantee that the Cerami sequences are bounded in the space $ X $, which norm is not invariant under translations. They have proved that it happens module translation. However, it makes difficult to apply multiplicity theorems since we do not have a Cerami condition being satisfied. In \cite{[cjj]}, the boundedness is obtained considering radial functions and prescribed norm, while in \cite{[6]} they have used subgroups of $ O(2) $. In this sense, in our work we provide a new technical result that allow us to get multiplicity for problem (\ref{P}), Lemma \ref{lkey}, in which we provide an argument that guarantee, up to a subsequence, the boundedness of Cerami sequences. 

We recall that a sequence $ (\un) \subset X $ is said to be a Palais-Smale sequence for the functional $ I $ in the level $ c $, $ (PS)_c $, if it satisfies 
$$ I(\un) \RA c \mbox{ \ \ and \ \ } ||I'(\un)||_{X'}\RA 0 \mbox{ , as \ } n \RA + \infty , $$ 
and $ (\un) $ is said to be a Cerami sequence for the functional $ I $ in the level $ c $, $ (PSC)_c $, if it verifies
$$ I(\un) \RA c \mbox{ \ \ and \ \ } ||I'(\un)||_{X'}(1+||\un||_X)\RA 0 \mbox{ , as \ } n \RA + \infty . $$
Clearly, any $ (PSC) $ sequence is a $ (PS) $ sequence.
  
In order to simplify some calculations and make the notation more concise, we will write
$$ R(\alpha , s) = \exp(\alpha |s|^{\frac{N}{N-1}}) - S_{N-2}(\alpha , s) = \sum\limits_{N-1}^{+\infty} \dfrac{\alpha^k}{k!} |s|^{\frac{N}{N-1}k} , $$
where $S_{N-2}(\alpha_0 , s) = \sum\limits_{k=0}^{N-2}\dfrac{\alpha_{0}^{k}}{k!}|s|^{\frac{N}{N-1}k}$. 

We recall that a funcion $ h $ is said to have \textit{subcritical} exponential growth at $ +\infty $, if
$$
\lim\limits_{s\RA + \infty}\dfrac{h(s)}{R(\alpha , s)} = 0 \textrm{ \ , for all \ } \alpha >0 ,
$$
and we say that $ h $ has $ \alpha_0 $-\textit{critical} exponential growth at $ +\infty $, if
$$
\lim\limits_{s\RA + \infty}\dfrac{h(s)}{R(\alpha , s)} = \left\{ \begin{array}{ll}
0, \ \ \ \forall \ \alpha > \alpha_{0} \\
+\infty , \ \ \ \forall \ \alpha < \alpha_0
\end{array} \right. .
$$

We list below the hypotheses assumed over $ g $ and $ f $, respectively. Then, we make a brief explanation why any of them are needed. The conditions over $ f $ are usual in works involving Moser-Trudinger inequalities, such as \cite{[doO], [zhang], [lam]}, and, in this case, are majorly based in \cite{[pucci], [boer3]}.

$$ g\in C(\R , \R), g(0)=0, g \mbox{ is odd and } \ G(s) \geq 0 \mbox{ \ , for all \ } s\in \R . \leqno{(g_1)} $$
$ (g_2) \mbox{ \ Given } \varepsilon > 0 , \mbox{ for some } \sigma \in (p, \sigma_0) , \mbox{ where } \sigma_0 = \frac{N}{p'}+1, \mbox{ there exists a constant } b_1 = b_1(\varepsilon , \sigma) > 0 \mbox{ such that }   $
$$|g(s)| \leq \varepsilon |s|^{p-1}+b_1|s|^{\sigma - 1} \ , \ \forall \ s \in \R .  $$ 
$$ \mbox{ There exists }\ \theta_1 \geq  \dfrac{N}{2} \ \mbox{ such that}\   g(s)s = \theta_1 G(s), \ \mbox{for all }\  s\in \R . \leqno{(g_3)}$$
$$
\mbox{There exist constants } b_2, b_4 > 0 \mbox{ such that } G(s) \geq b_2 |s|^N + b_4 |s|^p , \ \forall s \in \R . \leqno{(g_4)}
$$
$$ \mbox{ For some constant } \nu >0 , G(s-\tau) \leq \nu(g(s) - g(\tau))(s-\tau), \mbox{ for all } s, \tau \in \R . \leqno{(g_5)}$$ 

A prototype of function satisfying conditions $ (g_1)-(g_5) $ is given by the odd extension of $ g:[0, + \infty)\RA [0, + \infty) $ defined as $ g(s)=(\beta+1)|s|^\beta $, for $ s \geq 0 $, with $ \beta \in (p-1, \sigma_0 -1) $. For $ f $, we ask

$$ f\in C(\R , \R), f(0)=0, \mbox{has critical exponential growth and} \ F(t) \geq 0 \mbox{ \ , for all \ } t\in \R . \leqno{(f_1)} $$
$$ \lim\limits_{s\RA 0} \dfrac{f(s)}{|s|^{N-1}}=0.  \leqno{(f_2)}$$  
$$ \mbox{ there exists }\ \theta_2 \geq  2\theta_1 \ \mbox{ such that}\   f(s)s \geq \theta_2 F(s) > 0, \ \mbox{for all }\  s> 0. \leqno{(f_3)}$$
$$ \mbox{ there are}\  q>2N \ \mbox{ and}\  C_q> \dfrac{K_{10}^{\frac{N-1}{\beta_1}}}{\rho_{0}^{\frac{N}{\beta_1}}} \ \mbox{ such that}\  F(s) \geq C_q |s|^q , \ \mbox{for all}\  s>0 , \leqno{(f_4)}$$ 
where $ K_{10}=K_{10}(q, p, N)>0 $, $ \beta_1 $ and $ \rho_0>0 $ will be defined in Corollary \ref{l10}.

We start defining a functional $ J_t : \Lt \RA [0, + \infty] $ by
$$
J_t(u) = \intR \ln(1+|x|)|u(x)|^t dx .
$$
In the cases $ t=p $ and $ t=N $ we will denote, respectively, by $ J_p(u)=||u||_{\ast , p}^{p} $ and  $ J_N(u)=||u||_{\ast , N}^{N} $ because, as we will see, they will be norms in the spaces treated. Hence, inspired by \cite{[21]}, we define the spaces
$$
X_p = \{ u \in \Wump \ ; \ ||u||_{\ast, p} < + \infty\} \mbox{ \ \ and \ \ } X_N = \{ u \in \WumN \ ; \ ||u||_{\ast, N} < + \infty\},
$$
endowed with the norms $ ||\cdot||_{X_p}=||\cdot||_{1, p}^{p}+||\cdot||_{\ast, p}^{p} $ and  $||\cdot||_{X_N}=||\cdot||_{1, N}^{N}+||\cdot||_{\ast, N}^{N} $, respectively, where $ ||\cdot||_{1, t} $ stands for the usual norm of $ W^{1, t}(\RN) $, $ t>1 $. Set $ X=X_p \cap X_N $ with the norm $ ||\cdot||_X = ||\cdot||_{X_p}+||\cdot||_{X_N} $ and $ W= \Wump \cap \WumN $ with the norm $ ||\cdot||=||\cdot||_{1, p}+||\cdot||_{1, N} $. Sometimes, for simplicity, we will denote  $ ||\cdot||_\ast \coloneqq ||\cdot||_{\ast, p}^{p}+||\cdot||_{\ast, N}^{N} $.

The usual condition $ (g_3) $ is considered as an inequality instead of the equality. This is imposed in order to get the boundedness of Cerami sequences, Lemma \ref{l13}, in which appears the logarithmic term, that can changes sign. Moreover, we also need an Ambrosetti-Rabinowitz-type condition over $ f $, $ (f_3) $, to guarantee that the desired sequences are bounded. The arguments made in \cite{[6], [cjj], [10]} can not adapt immediately, because we need a bit more information about $ ||\nabla \cdot||_{N} $ in order to get the boundedness of the integral involving the exponential term.

Condition $ (g_4) $ is needed to obtain Lemma \ref{l9}, which makes possible to obtain boundedness and convergence inside the space $ X $. 

Finally, condition $ (g_5) $ makes possible to obtain the key Proposition to prove Theorem \ref{t1}, Proposition \ref{p2}, since the logarithmic changes sign and we need a way to control this term in order to get strong convergence in $ W $ and $ X $.  

Now we are ready to enunciate our first main result, which guarantees the existence of a ground state solution and a Mountain-Pass solution for (\ref{P}). 

\bt\label{t1}
Assume $ (f_1)-(f_4) $, $ (g_1)-(g_5) $, $ q>2N $ and that $ C_q>0 $ is sufficiently large. Then,
\begin{itemize}
\item[(i)] problem (\ref{P}) has a solution $ u\in X\setminus \{0\} $ verifying $ I(u)=c_{mp} $, where 
\begin{equation}\label{cmp}
c_{mp}= \inf\limits_{\gamma \in \Gamma} \max\limits_{t\in [0, 1]} I(\gamma(t)) ,
\end{equation}
with $ \Gamma = \{ \gamma\in C([0, 1], X) \ ; \ \gamma(0)=0 \ , \ I(\gamma(t))< 0 \} $.

\item[(ii)] Problem (\ref{P}) has a nontrivial ground state solution, i.e there exists $ u\in X\setminus \{0\} $ such that $ I(u)=c_g = \inf\{I(v) \ ; \ v\in X\setminus\{0\} \mbox{ \ is a solution of (\ref{P}}) \} $. 
\end{itemize}
\et

For the second main result, we are concerned with multiplicity of solutions. To do so, we will apply a version of the Symmetric Mountain-Pass Theorem, due to Rabinowtiz \cite{[albuquerque]}. For this purpose, we need the functional $ I $ to be even, so we replace the condition $ (f_1) $ by the condition below.
$$ f\in C(\R , \R), \ f \mbox{ \ is odd, has critical exponential growth and} \ F(s) \geq 0 \mbox{, for all \ } s\in \R . \leqno{(f_1 ')}$$
It is interesting to observe that the idea used in \cite{[boer3]} to obtain multiplicity cannot be applied to this, since we do not have the desired geometry, for instance, the homogeneity, even considering more growth conditions over $ f $ and $ g $.

Hence, we are able to enunciate the second main result. 
 
\bt\label{t2}
Suppose $ (f_1 ')-(f_4) $, $ (g_1)-(g_5) $, $ q> 2N $ and that $ C_q>0 $ is sufficiently large. Then, problem (\ref{P}) has infinitely many solutions. 
\et

Throughout this paper, we fix the values $ r_1, r_2 > 1 $, with $ r_1 \sim 1 $ and $ \frac{1}{r_1}+\frac{1}{r_2}=1 $, and we will use the following notations: $ \Lt $ denotes the usual Lebesgue space with norm $ ||\cdot ||_t $ \ ; \ $ X' $ denotes the dual space of $ X $ \ ; \ $ B_r(x) $ is the ball centred in $ x $ with radius $ r>0 $ \ ;  \ $ K_1, K_2, ... $ denotes constants that appears in important estimatives \ ; \ $ C, C_1, C_2, ... $ will denote different positive constants that appear inside the calculations and whose exact values are not essential \ ; \ $ \IC \IC $ stands for compact embeddings. 

The paper is organized as follows: in Section 2 we present the framework necessary to study problem (\ref{P}), the spaces involved and the regularity of the associated functional. Section 3 is devoted to get some technical results, concerning boundedness and convergence of sequences in the spaces of interest, and the geometry of the associated functional. In Section 4, we finish the proof of the first main result, while in Section 5 we prove the second main result. Finally, Section 6 consists in the discuss of related problems using the techniques applied in this work.

\section{Framework for problem (\ref{P})}

In this section, we will present some properties of the spaces $ W $ and $ X $ and provide some technical results that guarantee that the associated functional is well-defined and of class $ C^1 $ over $ X $.

Before we start providing the results, we call attention to the fact that, once $ \frac{N}{2}< p $, $ N < p^{\ast} $, where $ p^{\ast} = \frac{Np}{N-p} $ is the Sobolev critical exponent.

\bp\label{p1}
Concerning the space $ X $ we have the following embeddings:
\begin{itemize}
\item[(a)] $ X\IC \Wump $, $ X \IC \WumN $, $ X \IC W $ and $ X \IC \Lt $, for all $ t \geq p $.

\item[(b)] $X \IC \IC \Lt $, for all $ t \geq p $.
\end{itemize}
\ep
\begin{proof}
For item (a) one only needs to be reminded that $ \Wump \IC \Lt $, for all $ t\in [p, p^\ast ] $, and $ \WumN \IC \Lt $, for all $ t \geq N $. From this, the definition of $ ||\cdot|| $ and the fact that $ N < p^\ast $, the assertion follows. For item (b), we use the fact that $ \ln(1+|x|) $ is coercive, the interpolation inequality for $ \Lt $ spaces, the fact that $ N < p^\ast $ and a diagonal argument.
\end{proof}

We say that $ u\in X $ is a weak solution for (\ref{P}) if it verifies
\begin{align*}
& \intR [|\nabla u|^{p-2}\nabla u \nabla \varphi + |u|^{p-2}u\varphi ] dx + \intR  [|\nabla u|^{N-2}\nabla u \nabla \varphi + |u|^{N-2}u\varphi ] dx \\
& + \intRR \ln(|x-y|)G(u(x))g(u(y))\varphi(y) dx dy - \intR f(u)\varphi dx = 0 \ , \ \forall \ \varphi \in X .
\end{align*}
In this sense, we seek for an associated functional with (\ref{P}). First of all, inspired by \cite{[6]}, we define three auxiliary functionals, $V_1: W \RA [0, \infty],$ $V_2: Y \RA [0, \infty) $ and $V_0: W \RA \R \cup \{\infty\},$ given by 
\beqa
u \mapsto V_1(u)=\intRR \ln(1+|x-y|)G(u(x))G(u(y)) dx dy ,
\eeqa
\beqa
u \mapsto V_2(u, v)=\intRR \ln\left(1+\dfrac{1}{|x-y|}\right)G(u(x))G(u(y)) dx dy ,
\eeqa
\beqa
u \mapsto V_0(u, v)=V_1(u, v)-V_2(u, v)=\intRR \ln(|x-y|)G(u(x))G(u(y))dx dy. 
\eeqa
These definitions are understood to being over measurable function $u, v: \RN \RA \R $, such that the integrals are defined in the Lebesgue sense, where $ Y \coloneqq L^{\frac{2N}{2N-1}p}(\RN)\cap L^{\frac{2N}{2N-1}\sigma}(\RN)\cap L^{\frac{2N}{2N-1}(\sigma - 1)p'} (\RN)$, with $ \sigma \in (p, \sigma_0) $, where $ \sigma_0 = \frac{N}{p'}+1 $. Then, we consider the Euler-Lagrange functional associated with (\ref{P}), $ I: X \RA \R $, given by
$$
I(u) = \dfrac{1}{p}||u||_{1,p}^{p}+\dfrac{1}{N}||u||_{1, N}^{N}+ \dfrac{1}{2}V_0(u) - \intR F(u) dx .
$$
The aim of the following results in this section is to guarantee that $ I $ is well-defined for all $ u \in X $ and that $ I\in C^{1}(X, \R) $. We start focusing in the term involving the logarithm. 

\bl\label{l1}
Let $ r\in (p, N) $. Then, $ J_r(u) \leq ||u||_{\ast, p}^{p}+||u||_{\ast, N}^{N} < + \infty $, for all $ u\in X $.
\el
\begin{proof}
Consider the set $ \Omega_1=\{x\in \RN \ ; \ |u(x)|< 1 \} $ and $ \Omega_2=\{x\in \RN \ ; \ |u(x)|\geq 1 \}$. Since $ r \in (p, N) $, we have $ |u(x)|^r \leq |u(x)|^p $, for all $ x \in \Omega_1 $, and $ |u(x)|^r \leq |u(x)|^N $, for all $ x\in \Omega_2 $. Then, 
$$
J_r(u) \leq ||u||_{\ast, p}^{p}+||u||_{\ast, N}^{N} < + \infty \ , \ \forall \ u \in X .
$$
\end{proof}

\bo\label{obs1}
\begin{itemize}
\item[(1)] From condition $ (g_2) $, we obtain that $ \sigma \in (p, N) $ and $ (\sigma - 1)p' \in (p, N) $. Hence, Lemma \ref{l1} is applicable to these values. 

\item[(2)] Given $ \varepsilon >0 $, for some $ \sigma \in (p, \sigma_0) $, from conditions $ (g_2) $ and $ (g_4) $, there exists a constant $ b_3=b_3(\varepsilon, \sigma) >0 $  such that
\begin{equation}\label{eq1}
b_2 |s|^N + b_4 |s|^p \leq G(s) \leq \varepsilon |s|^p + b_3 |s|^\sigma \ , \ \forall \ s \in \R .
\end{equation}

\item[(3)]  Let $ u\in W $. Then, from Lemma \ref{l1} and $ (g_2) $, there exists a constant $ K_1=K_1(p, \sigma)>0 $, such that
\begin{equation}\label{eq3}
\intR \ln(1+|x|)|g(u)||u| dx \leq K_1(||u||_{\ast, p}^{p}+||u||_{\ast, p}||u||_{\ast, N}^{\sigma_0 -1}) \ , \ \forall \ u \in W. 
\end{equation}

\item[(4)] We also recall that, for $ r > 0 $, $ \ln(1+r) \leq r $ and
\begin{equation}\label{eq2}
\ln(1+|x-y|) \leq \ln(1+|x|)+\ln(1+|y|)  \ , \ \forall \ x, y \in \RN .
\end{equation} 
\end{itemize}
\eo

\bl\label{l2}
We have the following estimatives:
\begin{itemize}
\item[(a)] There exists a constant $ K_2=K_2(p, N, \sigma)>0 $ such that $$ V_2(u) \leq K_2 \left(||u||_{\frac{2N}{2N-1}p}^{2p}+ ||u||_{\frac{2N}{2N-1}\sigma}^{2\sigma}\right) \ , \ \forall \ u \in L^{\frac{2N}{2N-1}p}(\RN)\cap L^{\frac{2N}{2N-1}\sigma}(\RN) . $$

\item[(b)] There exists a constant $ K_3=K_3(p, N, \sigma)>0 $ such that
$$
V_1(u) \leq K_3 ( ||u||_{X}^{2p}+||u||_{X}^{p+N}+ ||u||_{X}^{p+\sigma}+ ||u||_{X}^{\sigma + N}) \ , \ \forall \ u \in X.
$$
\end{itemize} 
\el
\begin{proof}
For item (a), we use equation (\ref{eq1}), the Hardy-Littlewood-Sobolev inequality (HLS)(found in \cite{[15]}), with $ \alpha = \beta = 0 $, $ \lambda = 1 $ and a natural choice for $ q $ and $ t $, that is $ q=t=\frac{2N}{2N-1}$, and standard arguments. In the other hand, we have item (b) combining equations (\ref{eq1}) and (\ref{eq2}) and Lemma \ref{l1} with $ r=\sigma $. 
\end{proof}

\bo\label{obs3}
From Lemma \ref{l2}, the functionals $ V_1, V_2 $ and $ V_0 $ are well-defined over $ X $. Moreover, the functional $ V_2 $ is well-defined over $ L^{\frac{2N}{2N-1}p}(\RN)\cap L^{\frac{2N}{2N-1}\sigma}(\RN) $. 
\eo

The next two lemmas are technical results needed to prove that $ V_1 $ and $ V_2 $ are of $ C^1 $ classes. 

\bl\label{l3}
Let $ (\varphi_n) \subset X $ and $ \varphi \in X $. Then, 
\begin{itemize}
\item[(a)] if $ \varphi_n \RA \varphi $ in $ W $ or $ \varphi_n \CF \varphi $ in $ X $, then there exists a subsequence $ (\varphi_{n_k})\subset (\varphi_n) $ and a function $ h\in W $ such that  $ \varphi_{n_k}(x) \RA \varphi(x) $ a.e. in $ \RN $ and $ |\varphi_{n_k}(x)|\leq h(x) $, for all $ k\in \N $ and a.e. in $ \RN $.

\item[(b)] if $ \varphi_n \RA \varphi $ in $ X $, then there exists a subsequence $ (\varphi_{n_k})\subset (\varphi_n) $ and a function $ h\in X $ such that  $ \varphi_{n_k}(x) \RA \varphi(x) $ a.e. in $ \RN $ and $ |\varphi_{n_k}(x)|\leq h(x) $, for all $ k\in \N $ and a.e. in $ \RN $.
\end{itemize}
\el 
\begin{proof}
The proof is a consequence of \cite[Lemma 2.1]{[boer3]} and can be done similarly.
\end{proof}

\bl\label{l4}
Let $ \un \RA u $ in $ X $. Then, up to a subsequence, we have
$$
\intR \ln(1+|x|)|G(\un)-G(u)|dx \RA 0 \mbox{ \ \ and \ \ } \intR \ln(1+|x|)|g(\un)-g(u)|^{p'} dx \RA 0 . 
$$ 
\el
\begin{proof}
The proof follows from Lemma \ref{l3}, equation (\ref{eq1}), condition $ (g_2) $, the Dominated Convergence Theorem and by standard arguments.
\end{proof}

Then, we are ready to guarantee that $ V_1 $ and $ V_2 $ are of $ C^1 $ classes.In order to make the paper concise, we only sketchy the proof. Remember that $ Y \coloneqq L^{\frac{2N}{2N-1}p}(\RN)\cap L^{\frac{2N}{2N-1}\sigma}(\RN)\cap L^{\frac{2N}{2N-1}(\sigma - 1)p'} (\RN)$.

\bl\label{l5}
\begin{itemize}
\item[(a)] Let $ u, v \in X $. Then, the Gateuax derivative of $ V_1 $ is given by 
$$
V_1'(u)(v)=2 \intRR \ln(1+|x-y|)G(u(x))g(u(y))v(y) dx dy .
$$
Moreover, $ V_1 \in C^{1}(X, \R) $.

\item[(b)] Let $ u, v \in Y $. Then, the Gateuax derivative of $ V_2 $ is given by 
$$
V_2'(u)(v)=2 \intRR \ln\left(1+\dfrac{1}{|x-y|}\right)G(u(x))g(u(y))v(y) dx dy .
$$
Moreover, $ V_2 \in C^{1}(Y\cap L^{\frac{2N}{2N-1}}(\RN), \R) $.
\end{itemize}
\el
\begin{proof}
\textbf{(a)} We apply standard arguments, making use of the Mean Value Theorem, condition $ (g_2) $, equation (\ref{eq1}), Hölder inequality, Lemma \ref{l1}, the Dominated Convergence Theorem, Lemma \ref{l4} and equation (\ref{eq3}). Moreover, we get that
$$
||V_1'(u)||_{X'}\leq K_4 (||u||_{X}^{\sigma}+||u||_{X}^{2p-1}) \ , \ \forall \ u \in X ,
$$
where $ K_4=K_4(p, \sigma, N)>0 $.

\textbf{(b)} In light of Remark \ref{obs1}-(4), we use the Mean Value Theorem, the HLS inequality repeatedly, equation (\ref{eq1}), condition $ (g_2) $, the Dominated Convergence Theorem and standard arguments. Moreover, we obtain
$$
||V_2'(u)||_{Y'}\leq K_5 (||u||_{Y}^{3\sigma-1}+||u||_{Y}^{2p+\sigma-1}) \ , \ \forall \ u \in Y ,
$$
where $ K_5=K_5(p, \sigma, N)>0 $. The reader should observe that the intersection $Y\cap L^{\frac{2N}{2N-1}}(\RN)$ is needed only to guarantee that $ V_2' $ is continuous, since it will be necessary get some convergence related to $ L^{\frac{2N}{2N-1}}(\RN) $. 
\end{proof}

\bo\label{obs5}
\begin{itemize}
\item[(1)]From condition $ (g_3) $, $ V_2 '(u)(u)= \theta_1 V_2(u) $, for all $ u\in Y $. Hence, from Lemma \ref{l2}-(a), we conclude that there exists a constant $ K_6=\theta_1 K_2 >0 $ such that
$$
V_2 '(u)(u) \leq K_6\left(||u||_{\frac{2N}{2N-1}p}^{2p}+ ||u||_{\frac{2N}{2N-1}\sigma}^{2\sigma}\right) \ , \ \forall \ u \in Y .
$$

\item[(2)] From condition $ (g_3) $, $ V_1 '(u)(u)= \theta_1 V_1(u) $, for all $ u\in X $. Therefore, from Lemma \ref{l2}-(b), there exists a constant $ K_7 = \theta_1 K_3 >0 $ such that
$$
V_1 '(u)(u) \leq K_7 ( ||u||_{X}^{2p}+||u||_{X}^{p+N}||u||_{X}^{p+\sigma}+ ||u||_{X}^{\sigma + N}) \ , \ \forall \ u \in X .
$$
\end{itemize}
\eo

\bc\label{c1}
The functionals $ V_1, V_2 $ and $V_0$ belongs to $ C^1(X, \R) $. 
\ec

Next, we will focus on the term involving the nonlinearity $ f $. We start remembering the reader the famous Moser-Trudinger lemma, which makes possible to study such type of nonlinearities.

\bl\label{l6}
(Moser-Trudinger Lemma \cite{[doO], [17]}) Let $ N\geq 2 $, $ \alpha >0 $ and $ u\in \WumN $. Then, 
$$
\intR [ \exp(\alpha |u|^{N'}) - S_{N-2}(\alpha , u)] dx < \infty ,
$$
where $ S_{N-2}(\alpha, u) = \sum\limits_{k=0}^{N-2}\dfrac{\alpha^{k}}{k!}|u|^{N' k} $. Moreover, if $ ||\nabla u||_{N}^{N} \leq 1 $, $ ||u||_N \leq M < \infty $ and $ \alpha < \alpha_N = N\omega_{N-1}^{\frac{1}{N-1}} $, where $ \omega_{N-1} $ is the $ (N-1) $-dimensional measure of $ (N-1) $-sphere, then there exists a constant $ C_0 = C(\alpha , N , M) $ such that
$$
\intR [ \exp(\alpha |u|^{N'}) - S_{N-2}(\alpha , u)] dx \leq C(\alpha , N , M) = C_0 .
$$
\el

The next technical lemma allow us to conclude that $ R(\alpha , u)^l \in L^{1}(\RN) $, for $ \alpha > 0 , u \in X $, $ l\geq 1 $. 

\bl\label{l7}
(\cite[Lemma 2.3]{[w6]}) Let $ \alpha > 0 $ and $ r>1 $. Then, for every $ \beta > r $, there exists a constant $ C_\beta = C(\beta) > 0 $ such that
$$
(\exp(\alpha |t|^{p'}) - S_{N -2}(\alpha ,t))^r \leq C_\beta (\exp(\beta \alpha |t|^{p'} - S_{N -2}(\beta \alpha , t)) ,
$$
with $ \frac{1}{p}+\frac{1}{p'}=1 $. 
\el

\bc\label{c2}
Let $ \alpha>0 $. Then, $ R(\alpha , u)^l \in L^{1}(\RN) $, for all $ u\in \WumN $ and $ l \geq 1 $. In particular, it holds for all $ u\in X $.
\ec

\bc\label{c3}
Let $ u\in \WumN $, $ r>N $, $ l\geq 1 $, $ \beta > 0 $ and $ ||u||_{1, N}\leq M $, for $ M>0 $ sufficiently small. Then, there exists a constant $ K_8 = K_8(\beta , N, M, l) > 0 $ such that
$$
\intR |u|^r R(\beta , u)^l dx \leq K_8 ||u||_{t_0}^{r} ,
$$
for some $ t_0 > N $. Moreover, there exists a constant $ K_9 = K_9(\beta , N, M, l) > 0 $ such that
$$
\intR |u|^r R(\beta , u)^l dx \leq K_9 ||u||_{X}^{r} .
$$
\ec

The reader should observe that Corollary \ref{c3} asserts that is possible to obtain the estimates either controlling the norm $ ||\cdot||_{1, N} $, as we presented, or controlling the exponent $ \beta $ instead. Both forms are equally important inside of the work. In the sequence, we present two useful inequalities for the nonlinearity $ f $ and its primitive, both are valid for $ r>N $ although we will use it majorly for $ r=q>2N $.

\bo\label{obs4}
Given $ \varepsilon >0 $, $ r>N $ and $ \alpha > \alpha_0 $, from $ (f_2) $ and $ (f_1) $, for all $ u\in W $, there exists a constant $ b_5=b_5(\varepsilon , \alpha, r) >0 $ such that
\begin{equation}\label{eq4}
|f(u)| \leq \varepsilon |u|^{N-1} + b_5 |u|^{r-1}R(\alpha , u)
\end{equation}
Similarly, there exists a constant $ b_6=b_6(\varepsilon , \alpha, r) >0 $ satisfying
\beq\label{eq5}
|F(u)|\leq \varepsilon |u|^N + b_6 |u|^r R(\alpha ,u) .
\eeq
\eo

In order to see that the term involving $ F $ and $ f $ in the functional and in its derivative, respectively, are continuous, we need the following lemma. Its proof can be done similarly as in \cite[Lemma 2.5]{[boer3]}, so we omit it here.

\bl\label{l8}
Let $ (\un) \subset X $ and $ u\in X $ such that $ \un \RA u $ on $ W $. Then, we have
$$
\intR F(\un) \RA \intR F(u) \ \ , \ \ \intR f(\un)\un \RA \intR f(u)u \mbox{ \ \ and \ \ } \intR f(\un) v \RA \intR f(u)v \ , \ \forall \ v \in X .
$$
\el

Therefore, the results of this section guarantee that $ I $ is well-defined for all functions $ u\in X $ and that $ I\in C^{1}(X, \R) $. 

\section{Preliminary results and geometry of $I$}

We start this section providing some technical results concerning the boundedness and convergence of sequences in $ X $. Then, we verify that $ I $ has the mountain pass geometry and finish this section guaranteeing that certain sequence is bounded in $ W $ and under what conditions we can control the term involving $ f, F $, for this sequence.

\bl\label{l9}
Let $ u\in (\Lp \cap \LN ) \setminus \{0\} $ such that $ \un(x)\RA u(x) $ a.e. in $ \RN $ and $ (\vn)\subset \Lp \cap \LN $ bounded. Set 
$$
\omega_n = \intRR \ln(1+|x-y|)G(\un(x))G(\vn(y)) dx dy \ , \ \forall \ n \in \N .
$$
If $ \sup\limits_{n \in \N} \omega_n < +\infty $, then $ ||\vn||_{\ast} $ is bounded. Moreover, if $ \omega_n \RA 0 $ and $ ||\vn||_{\Lp \cap \LN}\RA 0 $, as $ n \RA + \infty $, then $ ||\vn||_{\ast}\RA 0 $.
\el 
\begin{proof}
First of all, from the Ergorov's Theorem, there exists $ R\in \N $, $ \delta > 0 $, $ n_0 \in \N $ and $ A\subset B_R $, such that $ A $ is measurable, $ \mu(A)>0 $ and $ \un(x) > \delta $, for all $ x\in A $ and for all $ n \geq n_0 $. From this and (\ref{eq1}), with $ \varepsilon = 1 $, together with $ \ln $ properties, we have
\begin{align*}
\omega_n & \geq \dfrac{b_{2}^{2}\delta^N \mu(A)}{2}\ds_{B_{2R}^{c}} \ln(1+|y|)|\vn(y)|^N dy + \dfrac{b_2 b_4 \delta^N \mu(A)}{2}\ds_{B_{2R}^{c}} \ln(1+|y|)|\vn(y)|^p dy \\
& + \dfrac{b_2 b_4\delta^p \mu(A)}{2}\ds_{B_{2R}^{c}} \ln(1+|y|)|\vn(y)|^N dy + \dfrac{b_{4}^{2}\delta^p \mu(A)}{2}\ds_{B_{2R}^{c}} \ln(1+|y|)|\vn(y)|^p dy \\
& \geq C_1 (||\vn||_{\ast , p}^{p} - \ln(1+2R)||\vn||_{p}^{p}+||\vn||_{\ast, N}^{N}- \ln(1+2R)||\vn||_{N}^{N}) \ , \ \forall \ n \geq n_0 ,
\end{align*}
where $ C_1 =\frac{1}{2} \min\left\{\frac{b_{2}^{2}\delta^N \mu(A)}{2} ,\frac{b_2 b_4 \delta^N \mu(A)}{2} , \frac{b_2 b_4\delta^p \mu(A)}{2}, \frac{b_{4}^{2}\delta^p \mu(A)}{2} \right\}  $. Thus, $ ||\vn||_\ast $ is bounded. Moreover, from the above inequality, if $ \omega_n \RA 0 $ and $ ||\vn||_{\Lp \cap \LN}\RA 0 $, we conclude that $ ||\vn||_\ast \RA 0 $.
\end{proof}

\bl\label{l10}
Let $ (\un)\subset X $ such that $ \un \rightharpoonup u $ in $ X $. Then, 
$$
\lim\limits_{n\RA + \infty} \intRR \ln(1+|x-y|)G(\un(x))g(u(y))(\un(y) - u(y)) dx dy = 0 .
$$
\el
\begin{proof}
We basically use condition $ (g_2) $ and equation (\ref{eq1}). The only term that appears in the calculations and we need to be careful is
$$
\intRR \ln(1+|y|) |\un(x)|^\sigma |u(y)|^{\sigma - 1}|\un(y)-u(y)| dxdy .
$$
But it can be dealt analogously as \cite[Lemma 3.1]{[boer3]}. 
\end{proof}

\bl\label{l11}
There exists $ \rho > 0 $ satisfying
\begin{equation}\label{eq6}
m_\beta = \inf \{ I(u) \ ; \ u\in X \ , \ ||u||=\beta \} \mbox{ , for all \ } \beta \in (0, \rho]
\end{equation}
and
\begin{equation}\label{eq7}
n_\beta = \inf \{ I'(u)(u) \ ; \ u\in X \ , \ ||u||=\beta \} \mbox{ , for all \ } \beta \in (0, \rho] .
\end{equation}
\el
\begin{proof}
Let $ u\in X \setminus \{0\} $, $ \alpha > \alpha_0 $ and $ r_1, r_2 >1 $, with $ r_1 \sim 1 $ and $ \frac{1}{r_1}+\frac{1}{r_2}=1 $, such that $ r_1 \alpha ||u||_{1, N}^{N'}< \alpha_N $. Then, from equation (\ref{eq5}), Lemma \ref{l2}-(a) and Lemma \ref{l6}, we have
\begin{align*}
I(u) & \geq \dfrac{1}{N}||u||_{1, p}^{N}+ \dfrac{1}{N}||u||_{1, N}^{N}- \dfrac{K_2}{2}||u||_{\frac{2N}{2N-1}p}^{2p}- \dfrac{K_2}{2}||u||_{\frac{2N}{2N-1}\sigma}^{2\sigma}-\varepsilon||u||_{N}^{N}-b_6 K_\alpha ||u||_{t_0}^{q} \\
& \geq \dfrac{||u||^N}{2^{N-1}N}[1-\varepsilon C_1 - C_2||u||^{2p-N}-C_3||u||^{2\sigma - N}-C_4 ||u||^{q-N}] .
\end{align*} 
Observe that, since $ \dfrac{N}{2}<p<\sigma $ and $ q>2N $, $ 2p-N > 0 $, $2\sigma - N>0$ and $ q-N >0 $. Thus, for $ \varepsilon, \rho>0 $ sufficiently small, we have $ I(u) > 0 $. On the other hand, from equation (\ref{eq5}), Remark \ref{obs5}-(1) and Lemma \ref{l6}, follows that
$$
I'(u)(u) \geq \dfrac{||u||^N}{2^{N-1}}[1-\varepsilon C_5 - C_6||u||^{2p-N}-C_7||u||^{2\sigma - N}-C_8 ||u||^{q-N}] .
$$
Hence, for $ \varepsilon, \rho > 0 $ sufficiently small, we conclude that $ I'(u)(u)>0 $.
\end{proof}

\bl\label{l12}
Let $ u\in X\setminus \{0\} $, $ t>0 $ and $ q> 2N $. Then,
$$
\lim\limits_{t\RA 0} I(tu) = 0 \ \ , \ \ \sup\limits_{t>0} I(tu) < +\infty \ \ \mbox{and} \ \ I(tu)\RA - \infty \ , \ \mbox{as} \ t \RA + \infty .
$$
\el
\begin{proof}
Let $ u\in X\setminus\{0\} $, $ t>0 $ and $ q>2N $. From $ (f_4) $ and Lemma \ref{l2}-(b), we have
\begin{align*}
I(tu) & \leq \dfrac{t^p}{p}||u||_{1, p}^{p} + \dfrac{t^N}{N}||u||_{1, N}^{N}+ K_3 t^{2p}||u||_{X}^{2p}+ K_3 t^{p+N}||u||_{X}^{p+N}+K_3 t^{p+\sigma}||u||_{X}^{p+ \sigma} \\
& + K_3 t^{\sigma + N}||u||_{X}^{\sigma +N}- C_q t^q ||u||_{q}^{q}\RA - \infty \ , \mbox{ as } t \RA + \infty .
\end{align*}
Now, consider $ t>0 $ such that $ r_1 \alpha t^{N'} ||u||_{1, N}^{N'}< \alpha_N  $. Then, from (\ref{eq5}) and Lemma \ref{l6}, we obtain that $ \intR F(tu) dx \RA 0 $. From Lemma \ref{l2}, $ V_0(tu) \RA 0 $ as well. Therefore, $ I(tu) \RA 0 $ as $ t\RA 0 $. Finally, from the fact that $ I $ is $ C^1 $, we conclude that $ \sup\limits_{t>0} I(tu) < + \infty $.
\end{proof}

\bo\label{obs6}
From Lemmas \ref{l11} and \ref{l12}, we conclude that $ I $ has the mountain pass geometry. Thus, the value $ c_{mp} $, defined in (\ref{cmp}), is well-defined and satisfies $ 0< m_\rho \leq c_{mp} <+\infty $. Moreover, there exists a Cerami sequence for $ I $ in level $ c_{mp} $.
\eo

Consider a sequence $ (\un)\subset X $ satisfying
\begin{equation}\label{eq8}
\exists \ d > 0 \mbox{ \ s.t. \ } I(\un) \leq d \ , \ \forall \ n \in \N \mbox{ \ and \ } ||I'(\un)||_{X'}(1+||\un||_X) \RA 0 \ , \ \mbox{ as \ } n \RA + \infty .
\end{equation}

\bl\label{l13}
Let $ (\un) \subset X $ satisfying (\ref{eq8}). Then, $ (\un) $ is bounded in $ W $.
\el 
\begin{proof}
From (\ref{eq8}), $ (f_3) $ and $ (g_3) $, we have
$$
d + o(1) \geq I(\un) - \dfrac{1}{2\theta_1}I'(\un)(\un) \geq \left(\dfrac{1}{p}- \dfrac{1}{2\theta_1}\right)||\un||_{1, p}^{p}+ \left(\dfrac{1}{N}- \dfrac{1}{2\theta_1}\right)||\un||_{1, N}^{N} \ , \ \forall \ n \in \N .
$$
Then, 
$$
\left[\left(\dfrac{1}{N}- \dfrac{1}{2\theta_1}\right)d\right]^{\frac{1}{N}}+ \left[\left(\dfrac{1}{p}- \dfrac{1}{2\theta_1}\right)d\right]^{\frac{1}{p}}+o(1) \geq ||\un|| \ , \ \forall \ n \in \N .
$$
Therefore, $ (\un) $ is bounded in $ W $.
\end{proof}

\bo\label{obs7}
From Lemma \ref{l13}, for a sequence satisfying (\ref{eq8}), with $ d=c_{mp} $, or a Cerami sequence in the level $ c_{mp} $, we conclude that $$ ||\un||_{1, N}^{N'} \leq o(1) + \left(\dfrac{1}{N}- \dfrac{1}{2\theta_1}\right)^{\frac{N'}{N}}c_{mp}^{\frac{N'}{N}} \ , \ \forall \ n \in \N . $$
Thus, up to a subsequence, we obtain $ ||\un||_{1, N}^{N'} \leq \left(\frac{1}{N}\right)^{\frac{1}{N-1}}c_{mp}^{\frac{1}{N-1}} $, for all $ n \in \N $.
\eo

\bl\label{l14}
There exists a constant $ K_{10} = K_{10} (q, p, N)>0 $ such that $ c_{mp} \leq \frac{K_{10}}{C_{q}^{\beta^1}} $, where $ \beta_1 = \frac{p}{q-p} $
\el
\begin{proof}
First of all, define the set $ \A = \{u \in X \ ; \ u \neq 0 \mbox{ \ and \ } V_0(u) \leq 0\} $. We need to verify that $ \A \neq \emptyset $. In order to do so, we define the set $ \Omega = \{ (x, y) \in \RN\times \RN \ ; \ |x-y| \geq 1 \} \subset \RN\times \RN $ and, for $ t > 0 $, define $ u_t(x) = t^2 u(tx) $, for all $ x \in \RN $. Moreover, to make the notation concise, set
$$
\varphi_{s, r}(u) \coloneqq \iint\limits_{\Omega}\ln(|x-y|)|u(x)|^{s}|u(y)|^r dx dy \ , \ \forall \ u \in X.
$$
Then, from (\ref{eq1}) and change of variables, we have
\begin{align*}
V_0(u) & \leq \iint\limits_{\Omega}\ln(|x-y|)G(u(x))G(u(y)) dx dy \\
& \leq [1-\ln t]t^{4p-2}\varphi_{p, p}(u) +  2b_3 [1-\ln t]t^{2p+2\sigma - 2}\varphi_{p, \sigma}(u) + b_{3}^{2}[1-\ln t] t^{4\sigma -2}\varphi_{\sigma, \sigma}(u)\RA - \infty ,  
\end{align*}
as $ t \RA + \infty $. Hence, $ \A \neq \emptyset $. Now, from \cite[Lemma 2.1]{[pucci]}, there exists a constant $ C=C(p, q, N)>0 $ such that $ ||u|| \geq C ||u||_{q} $. Thus, for any $ v \neq 0 $, makes sense to define
$$
S_q (v) = \dfrac{||v||}{||v||_q} \geq C \mbox{ \ \ and \ \ } S_q = \inf\limits_{v \in \A} S_q(v) \geq \inf\limits_{v \in W \setminus \{0\}} S_q(v) > 0 .
$$
From Lemma \ref{l12}, for $ v\in \A $ and $ T>0 $ sufficiently large, we have $ I(Tv)<0 $. Consider $ \gamma: [0, 1]\RA X $ given by $ \gamma(t)=tTv $. Then, $ \gamma \in \Gamma $ and $ c_{mp} \leq \max\limits_{t \geq 0}I(tv) $.
Hence, for $ \psi \in \A $, we obtain
\begin{align*}
c_{mp} & \leq \max\limits_{t \geq 0}\left\{\dfrac{S_q(\psi)^p}{p}t^p ||\psi||_{q}^{p}+ \dfrac{S_q(\psi)^N}{p}t^N ||\psi||_{q}^{N} - C_q t^q ||\psi||_{q}^{q} \right\} \\
& = \max\limits_{s \geq 0}\left\{\dfrac{S_q(\psi)^p}{p}s^{p}+ \dfrac{S_q(\psi)^N}{p}s^{N} - C_q s^{q} \right\} \\
& \leq \max\limits_{s \geq 0} h_1(s) + \max\limits_{s \geq 0} h_2(s),  
\end{align*}
where $ h_1, h_2: [0, + \infty) \RA \R $ are given, respectively, by
$$
h_1(s) = \dfrac{S_q(\psi)^p+S_q(\psi)^N}{p} s^{p} - C_q s^{q} \mbox{ \ \ and \ \ } h_2(s) = \dfrac{S_q(\psi)^p+S_q(\psi)^N}{p}  s^{N} - C_q s^{q}
$$
Through some calculations and taking the infimum over all $ \psi \in \A $, one can find two constants $ C_1=C_1(p, q, N) >0 $ and  $ C_2=C_2(p, q, N) >0 $ such that $ \max\limits_{s \geq 0} h_1(s)\leq \frac{C_1}{C_{q}^{\beta_1}} $ and $ \max\limits_{s \geq 0} h_2(s) \leq \frac{C_2}{C_{q}^{\beta_1}} $. Therefore, considering $ K_{10} = 2\max\{C_1, C_2\} $, we get the result.
\end{proof}

The next corollary is an immediate consequence of Lemma \ref{l14} and Remark \ref{obs7}. It will allow us to apply the Moser-Trudinger Lemma for sequences under those hypotheses, once we can make $ \rho_0 $ as small as we need.

\bc\label{c4}
Let $ (\un) \subset X $ satisfying (\ref{eq8}), for $ d \in (0, c_{mp}] $, or being a Cerami sequence for $ I $ at level $ c_{mp} $, $ q> 2N $ and $ C_q>0 $ sufficiently large. Then, we can find $ \rho_0 $ sufficiently small, such that
$$
\limsup\limits_{n}||\un||_{1, N}^{N'} \leq \rho_{0}^{N'} .
$$
\ec

\bl\label{l15}
Let $ (\un)\subset X $ be bounded in $ W $ such that 
\begin{equation}\label{eq9}
\liminf \sup\limits_{y\in \mathbb{Z}^N} \ds_{B_2(y)}|\un(x)|^p dx > 0 .
\end{equation}
Then, there exists $ u\in W \setminus\{0\} $ and $ (\yn)\subset \mathbb{Z}^N $ such that, up to a subsequence, $ \yn \ast \un = \until \CF u $ in $ W $. 
\el
\begin{proof}
From (\ref{eq9}) and $ \liminf $ properties, there exists a constant $ C_1>0 $ such that, passing to a subsequence if necessary, 
$$
\sup\limits_{y\in \mathbb{Z}^N} \ds_{B_2(y)}|\un(x)|^p dx > C_1 \ , \ \forall \ n \in \N .
$$ 
Thus, for each $ n\in \N $, from $ \sup $ definition, there exists a sequence $ (y_{k}^{n})\subset \mathbb{Z}^N $ satisfying
$$
\lim\limits_{k\RA + \infty} \ds_{B_2(y_{k}^{n})}|\un(x)|^p dx = \sup\limits_{y\in \mathbb{Z}^N} \ds_{B_2(y)}|\un(x)|^p dx > C_1 .
$$
Then, for each $ n\in \N $, there exists $ k_{0}^{n}\in \N $ verifying
$$
\ds_{B_2(y_{k_{0}^{n}}^{n})}|\un(x)|^p dx > C_1 .
$$
Moreover, since $ (\un) $ is bounded in $ W $, there exists $ C_2>0 $ such that 
\begin{equation*}
C_2 > ||\un||^p \geq ||\un||_{p}^{p} \geq \ds_{B_2(y_{k_{0}^{n}}^{n})}|\un(x)|^p dx > C_1 .
\end{equation*}
Hence, we obtain a sequence indexed in $ n\in \N $ satisfying
\begin{equation}\label{eq10}
\left(\ds_{B_2(y_{k_{0}^{n}}^{n})}|\un(x)|^p dx\right) \subset [C_2 , C_1].
\end{equation}
Consider $ (\un) $ and $ (\yn) $ the subsequences obtained by the above construction. Let $ \until = \yn \ast \un $. So, as $ ||\cdot|| $ is $ \mathbb{Z}^N $-invariant, $ (\until) $ is bounded in $ W $. Therefore, there exists $ u\in W $ such that $ \until \CF u $ in $ W $. 

As a consequence, without loss of generality, we can assume that $ \until (x) \RA u(x) $ a.e. in $ \RN $. Let us show that $ u\neq 0 $ in $ W $. 

To begin with, since $ \Wump \IC W $ and the restriction operator from $ \Wump $ to $ W^{1, p}(B_2) $ is continuous, we have that  $ \until \big\vert_{B_2} \CF u\big\vert_{B_2} $ in $ W^{1, p}(B_2) $. Then, from Rellich-Kondrakov, $ \until \big\vert_{B_2} \RA u\big\vert_{B_2} $ in $ L^{p}(B_2) $. Thus, 
$$
C_1 < \ds_{B_2}|\until|^p dx = \ds_{B_2}|\until\big\vert_{B_2}|^p dx \RA \ds_{B_2}|u\big\vert_{B_2}|^p dx \leq ||u||_{p}^{p}.
$$
Hence, $ ||u||^p > C_1 $ which implies $ u\neq 0 $ in $ W $.
\end{proof}

\bl\label{l16}
Let $ (\un)\subset X $ be sequence satisfying either (\ref{eq8}), with $ d\in (0, c_{mp}] $, or being a Cerami sequence for $ I $ at level $ c_{mp} $, and that does not verify $ ||\un||\RA 0 $ and $ I(\un) \RA 0 $, $ q> 2N $. Then, 
$$
\liminf \sup\limits_{y\in \mathbb{Z}^N} \ds_{B_2(y)}|\un(x)|^p dx > 0 .
$$
\el
\begin{proof}
Suppose, by contradiction, that 
$$
\liminf \sup\limits_{y\in \mathbb{Z}^N} \ds_{B_2(y)}|\un(x)|^p dx = 0 .
$$
Then, by Lion's Lemma \cite{[lions]}, $ \un \RA 0 $ in $ \Lt $, for all $ t \in (p, p^\ast) $. In particular, for $ t=N $. Then, 
$$
\liminf \sup\limits_{y\in \mathbb{Z}^N} \ds_{B_2(y)}|\un(x)|^N dx = 0 
$$
and, once more by Lion's Lemma, we obtain $ \un \RA 0 $ in $ \Lt $, for all $ t \geq N $. Thus, from (\ref{eq4}) and Remark \ref{obs5}-(1), 
$$
||\un||_{1, p}^{p}+||\un||_{1, N}^{N}+V_1'(\un)(\un) = I'(\un)(\un) + V_2 ' (\un)(\un) + \intR f(\un) \un dx \leq o(1) + \varepsilon C . 
$$
Hence, $ ||\un||_{1, p}^{p}+||\un||_{1, N}^{N}+V_1'(\un)(\un)\RA 0 $, as $ n \RA + \infty $, $ \varepsilon \RA 0 $. Consequently, $ ||\un||\RA 0 $ and, from Remark \ref{obs5}-(2), $ V_1(\un) \RA 0 $. 

Therefore, from the continuous embedding and (\ref{eq5}), we conclude that $ I(\un) \RA 0 $, which contradicts the hypotheses, finishing the proof.
\end{proof}

\bc\label{c5}
Let $ (\un)\subset X $ under the hypotheses given in Lemma \ref{l16}. Then, $ (\until) \subset X $ is bounded.
\ec
\begin{proof}
From Lemma \ref{l16}, $\liminf \sup\limits_{y\in \mathbb{Z}^N} \ds_{B_2(y)}|\un(x)|^p dx > 0$. Then, from Lemma \ref{l15}, there exists a sequence of points $ (\yn) \subset \mathbb{Z}^N $ such that, up to a subsequence, $ \until \CF u \in W \setminus \{0\} $. From Lemma \ref{l13} and the invariance of $ ||\cdot|| $, $ (\until) $ is bounded in $ W $ and, by the continuous embeddings, $ (\until) $ is bounded in $ \Lt $, for all $ t \geq p $. Moreover, we can assume, without loss of generality, that $ \until(x)\RA u(x) $ a.e. in $ \RN $. 

Consequently, from the invariance of $ V_1 , V_1 ' $, (\ref{eq8}), Remark \ref{obs5}, (\ref{eq4}) and Corollary \ref{c3},
$$
\theta_1 V_1(\until) = V_1'(\un)(\un) \leq I'(\un)(\un) + V_2 ' (\un)(\un) +  \intR f(\un) \un dx \leq C \ , \ \forall \ n \in \N .
$$
Therefore, $ \sup\limits_{n \in \N} V_1(\until) < + \infty $ and, from Lemma \ref{l9}, we conclude that $ (\until) $ is bounded in $ X $.
\end{proof}

The next technical lemma is the key to obtain the multiplicity of solutions, Theorem \ref{t2}, since it makes possible to verify the condition $ (PSC) $ in some suitable levels. It is the first work dealing with logarithmic Choquard equations presenting a boundedness result for a subsequence of the original sequence, $ (\un) $, in the space $ X $.

\bl\label{lkey}
Let $ (\un)\subset X $ under the hypotheses given in Lemma \ref{l16}. Then, up to a subsequence, $ (\un) $ is bounded in $ X $.
\el
\begin{proof}
First of all, from Lemma \ref{l13} $ (\un) $ is bounded in $ W $ and from Lemma \ref{l15}, passing to a subsequence if necessary,there exists $ (\yn)\subset \mathbb{Z}^N $ such that, up to a subsequence, $ \until \CF u $ in $ W\setminus\{0\} $ and $ \until(x) \RA u(x) $ pointwise a.e. in $ \RN $. Moreover, from $ W\IC \Wump $ and $ W\IC \WumN $, $ \until \CF u $ in $ \Wump $ and in $ \WumN $ and, from the proof of Lemma \ref{l15}, one can see that $ u\neq 0 $ in both spaces $ \Lp $ and $ \LN $. We will make the proof only for $ p $, since for $ N $ will be very similar. Lets split the proof into two cases.

\noindent \textbf{Case 1:} $|\yn|\RA +\infty$.

Observe that, if $ |\yn|\RA +\infty $, $ |x+\yn|\RA + \infty $, for all $ x\in \RN $, and
$$
\left| \dfrac{1+|\yn|}{1+|x+\yn|}-1\right| = \dfrac{||\yn|-|x+\yn||}{1+|x+\yn|}\leq \dfrac{|x+\yn - \yn|}{1+|x=\yn|} = \dfrac{|x|}{1+|x+\yn|}\RA 0 .
$$
Then, for $ x\in \RN $,
$$
\ln(1+|\yn|)-\ln(1+|x+\yn|)=\ln \left(\dfrac{1+|\yn|}{1+|x+\yn|}\right) \RA 0 , n\RA + \infty .
$$
So, for $ \varepsilon > 0 $ sufficiently small, once $ \ln(1+|x+\yn|)\RA + \infty $ as $ n \RA + \infty $, there exists $ n_0 \in \N $ and $ C> $ such that $ \varepsilon \leq C \ln(1+|x+\yn|) $ and 
$$
\ln(1+|\yn|)-\ln(1+|x+\yn|) \leq \varepsilon \ \Rightarrow \ \ln(1+|\yn|) \leq \varepsilon + \ln(1+|x+\yn|) \leq (1+C) \ln(1+|x+\yn|) , 
$$
for all $ n \geq n_0 $. Then, passing to a subsequence if necessary, there exists a constant $ C_1 > 0 $ such that $ \ln(1+|x+\yn|)\geq C_1 \ln(1+|\yn|) $, and
$$
||\until||_{\ast, p}^{p}=\intR \ln(1+|x|)|\un (x-\yn)|^p dx = \intR \ln(1+|x+\yn|) |\un(x)|^p dx \geq C_1 ||\un||_{p}^{p}\ln(1+|\yn|),
$$
for all $ n \in \N $. Now, since $ \until \RA u $ in $ \Lp \setminus \{0\} $ and $ ||\cdot||_p $ is $ \mathbb{Z}^N$-invariant, if $ ||\un||_p \RA ||u||_p >0 $. Hence, considering a subsequence if necessary, $ ||\un||_{p}^{p} \geq C_2 > 0 $, for all $ n \in \N $ and we obtain that 
$$
||\until||_{\ast, p}^{p} \geq C_2 \ln(1+|\yn|) , \forall \ n \in \N .
$$

\noindent \textbf{Case 2:} Suppose that $ (\yn)\subset \mathbb{Z}^N $ converges to $ y_0 \in \mathbb{Z}^N $. 

If $ y_0 =0 $, $ \until = \un $, nothing remains to be proved. So, we suppose that $ y_0 \neq 0 $. Then, up to a subsequence, $ \yn \equiv y_0 $. Observe that, given any $ \varphi \in C_{0}^{\infty}(\RN) $, $ y_0\ast \varphi \in C_{0}^{\infty}(\RN) $. For change of variables and $ \yn \equiv y_0 $, for each $ n\in \N $, 
$$
\intR \un(x) \varphi(x) dx = \intR \un(x-\yn) \varphi(x-\yn) dx = \intR \until(x)\varphi(x-y_0) dx = \intR \until (x)(y_0 \ast \varphi)(x) dx
$$
and
$$
\intR u(x) (y_0 \ast \varphi)(x) dx = \intR u(x)\varphi(x-y_0)= \intR u(x+y_0) \varphi(x) dx = \intR ((-y_0)\ast u) \varphi dx 
$$
Thus, since $ \until \CF u $ in $ \Wump \setminus\{0\} $,
$$
\intR \un \varphi dx = \intR \until (y_0 \ast \varphi) dx \RA \intR u (y_0 \ast \varphi) dx = \intR ((-y_0)\ast u) \varphi dx .
$$
Similarly, we obtain
$$
\intR |\nabla \un |^{p-2}\nabla \un \nabla \varphi dx \RA \intR |\nabla ((-y_0)\ast u)|^{p-2}\nabla ((-y_0)\ast u) \nabla \varphi dx .
$$
Hence, $ \un \CF ((-y_0)\ast u) \coloneqq u_0 $ and, since $ u \neq 0 $ in $ \Wump $, $ u_0 \neq 0 $ in $ \Wump $. Moreover, we can assume, without loss of generality, that $ \un(x) \RA u_0(x) $ pointwise a.e. in $ \RN $. 

Next, we will construct a suitable subset of $ \RN $.

Consider $ r $  the line passing through the origin and $ y_0 $. Then, define $ \Omega_0 $ as the open connected region between $ r $ and one of the axis, such that the angle between $ r $ and the axis is $ \leq \frac{\pi}{2} $. Since $ \Omega_0 $ is open, we can fix a point $ x_1 \in \Omega_0 $ and consider the domain $ \Omega = B_{\delta}(x_1)\subset \Omega_0 $, for $ \delta >0 $. 

\noindent \textbf{Claim:} There exists a constant $ C_3>0 $ and $ n_0\in \N $ such that $ ||\un||_{p, \Omega}^{p}\geq C_3 > 0 $, for all $ n \geq n_0 $.

Once $ X \hookrightarrow \hookrightarrow \Lp $ and $ \until \CF u $ in $ X $, $ \until \RA u $ in $ \Lp $.  Consequently, $ \until \RA u $ in $ L^p(\Omega) $ which implies $ ||\until||_{p, \Omega} \RA ||u||_{p, \Omega} $. Since $ ||\cdot||_p $ is $ \mathbb{Z}^N $ invariant we have that $ ||\un||_{p, \Omega} \RA ||u_0||_{p, \Omega} $ and, from a classical result, $ \un \RA u_0 $ in $ L^p(\Omega) $. 

Therefore, since $ u_0 \in \Lp \setminus \{0\} $, $ ||u_0||_{p, \Omega} >0 $ and there exists a constant $ C_3 > 0 $ and $ n_0 \in \N $ such that $ ||\un||_{p, \Omega}^{p} \geq C_3 > 0 $, for all $ n \geq n_0 $, proving the claim.

Observe that, for $ x\in \Omega $, we have $ |x+y_0| > |y_0| $, once the size of the diagonal of the parallelogram will be greater than its side $ |y_0| $.  

Hence, from Claim 1, $ \yn \equiv y_0 $ and $ \ln(1+|x+y_0|) \geq \ln(1+|y_0|) $, for all $ x\in \Omega $, we have
$$
||\until||_{\ast, p}^{p}  = \intR \ln(1+|x+\yn|)|\un(x)|^p dx \geq \ln(1+|y_0|)||\un||_{p, \Omega}^{p} \geq C_3 \ln(1+|\yn|) , \forall \ n \geq n_0 .
$$

From cases 1 and 2, we conclude that, up to a subsequence, there exists a constant $ C_4 >0 $ such that
\begin{equation}\label{eqk1}
||\until||_{\ast, p}^{p}= \intR \ln(1+|x+\yn|)|\un (x)|^p dx \geq C_4 \ln(1+ |\yn|) \ , \forall \ n \in \N .
\end{equation}
On the other hand, once $ \ln(1+|x+\yn|) \leq \ln(1+|x|)+\ln(1+|\yn|) $, for all $ n \in \N $, we have
\begin{equation}\label{eqk2}
||\un||_{\ast, p}^{p} = \intR \ln (1+|x+\yn|) \until^2(x) dx \leq ||\until||_{\ast, p}^{p} + \ln (1+|\yn|) ||\until||_{p}^{p} .
\end{equation}
From (\ref{eqk1}) and the fact that $ (\until) $ is bounded in $ X $, we have
$$
\ln (1+|\yn|) ||\until||_{p}^{p} \leq C_{4}^{-1}||\until||_{\ast, p}^{p}||\until||_{p}^{p} \leq C_5 ,
$$
for some constant $ C_5 > 0 $. Therefore, returning to (\ref{eqk2}), we obtain $ ||\un||_{\ast, p}^{p}\leq C_6 $, for some constant $ C_6>0 $. Analogously, there exists a constant $ C_7 > 0 $ such that $ ||\un||_{\ast, N}^{N}\leq C_7 $, for all $ n \in \N $.

Hence, since $ (\un) $ is already bounded in $ W $, we conclude that $ (\un) $ is bounded in $ X $.
\end{proof}

\section{Proof of Theorems \ref{t1}}

In this section we present the proof of the existence theorem. We begin proving a key proposition that give us the conditions to existence of nontrivial critical points for $ I $ is $ X $. 

\bp\label{p2}
Let $ q \geq 2N $ and $ (\un)\subset X $ a sequence either satisfying (\ref{eq8}), with $ d\in (0, c_{mp}] $ or being a Cerami sequence for $ I $ at level $ c_{mp} $. Then, passing to a subsequence, if necessary, only one between the alternatives below occurs:

\noindent \textbf{(a)} $ ||\un||\RA 0 $ and $ I(\un)\RA 0 $.

\noindent \textbf{(b)} There exists a function $ u\in X\setminus\{0\} $ such that $ \un \RA u $ in $ X $, for a non-trivial critical point $ u \in X $ of $ I $.
\ep
\begin{proof}
Suppose that item (a) does not happen. Then, from Lemma \ref{l13}, $ (\un) $ is bounded in $ W $. Also, from Lemmas \ref{l15} and \ref{l16} and Corollary \ref{c5}, there exist a sequence $ (\yn) \subset \mathbb{Z}^N $ and $ u_0 \in W \setminus\{0\} $ such that, up to a subsequence, $ \until \CF u_0 $ in $ W \setminus\{0\} $, $ (\until) $ is bounded in $ X $ and $ \until \CF u_0 $ in $ X \setminus\{0\} $. Moreover, from Proposition \ref{p1}, $ \until \RA u_0 $ in $ \Lt $, for all $ t \geq p $. 

Hence, from Lemma \ref{lkey}, up to a subsequence, $ (\un) $ is bounded in $ X $ and $ \un \CF u $ in $ X$. Once more, from Proposition \ref{p1},  $ \un \RA u $ in $ \Lt $, for all $ t \geq p $. If $ u=0 $ in $ X $, then $ u=0 $ in $ \Lp $ and from the $ \mathbb{Z}^N$-invariance of $ ||\cdot||_p $, we obtain that $ ||u_0||_p = 0 $, which is a contradiction. Then, $ u\in X\setminus\{0\} $.

Now, observe that, 
$$ |I'(\un)(\un - u)|\leq ||I'(\un)||_{X'}||\un - u||_X \leq C_1 ||I'(\un)||_{X'} \RA 0 \mbox{, as \ }  n \RA + \infty .$$

One can easily verify, similarly as \cite[Proposition 4.1]{[boer3]}, that
$$
\intR f(\un)(\un - u) dx \RA 0 \ , \ V_2 '(\un)(\un - u) \RA 0 \mbox{ \ \ and \ \ } \intR |\un|^{\beta-2}\un (\un - u) dx \RA 0 \ , 
$$
as $ n \RA + \infty $, with $ \beta = p, N $. Moreover, from condition $ (g_4) $, for all $ n \in N $, we get
\begin{align*}
V_1 '(\un)(\un- u) & = \intRR \ln(1+|x-y|)G(\un(x))g(\un(y))(\un(y)-u(y)) dx dy \\
& \geq \dfrac{1}{\nu}\intRR \ln(1+|x-y|)G(\un(x))G(\un(y)-u(y)) dx dy \\
& + \intRR \ln(1+|x-y|)G(\un(x))g(u(y))(\un(y)-u(y)) dx dy = A + B .
\end{align*}
From $ (g_1) $, $ A \geq 0 $ and, from Lemma \ref{l10}, $ B \RA 0 $, as $ n \RA + \infty $. Hence, $ V_1 '(\un)(\un - u) \geq o(1) $. As a consequence, 
\begin{small}
\begin{align*}
o(1) & = I'(\un)(\un - u) \\
& \geq \intR |\nabla \un |^{p-2}\nabla \un \nabla (\un - u) dx + \intR |\nabla \un |^{N-2}\nabla \un \nabla (\un - u) dx + V_1 '(\un)(\un - u) + o(1) \geq o(1) .
\end{align*}
\end{small}
Therefore, one can see that $ \un \RA u $ in $ W $ and $ V_1 '(\un)(\un - u) \RA 0 $. From Remark \ref{obs5}-(2) and Lemma \ref{l9}, $ \un \RA u $ in $ X\setminus\{0\} $.

Finally, it remains to prove that $ u $ is a critical point of $ I $. Let $ v \in X $. Then,
$$
|I'(u)(v)| = \lim |I'(\un)(v)| \leq ||v||\lim ||I'(\un)||_{X'} = 0 .
$$
\end{proof}

\begin{proof}[Proof of Theorem \ref{t1}]
Item \textbf{(a)} follows from Lemma \ref{l11} and Proposition \ref{p2}. For item \textbf{(b)}, we define the set $ \K =\{v \in X\setminus \{0\} \ ; \ I'(v)=0\} $, which is non-empty, since $ u $ obtained in item (a) belongs to $ \K $. Thus, we can consider a sequence $ (\un)\subset \K $ satisfying $ I(\un)\RA c_g = \inf\limits_{v\in K} I'(v) $.

Now, observe that $ c_g \in [-\infty , c_{mp}] $. If $ c_g = c_{mp} $ the proof is done. Otherwise, if $ c_g < c_{mp} $, from  the definition of $ \K $ and Proposition \ref{p2} there exists a function $ u_1 $ in $ X $ that is a non-trivial critical point of $ I $ in $ X $. Then, we conclude that $ u\in \K $ and 
$
I(u)= c_g .
$
Particularly, we see that $ c_g > -\infty $.
\end{proof}

\section{Proof of Theorem \ref{t2}}

This section is devoted to prove Theorem \ref{t2}. To do so, we will apply a version of Symmetric Mountain-Pass Theorem, due to Rabinowitz \cite{[ambrosseti]} (see also \cite{[bartolo], [silva]}).

\bt\label{t3}
(\cite[Theorem 4.1]{[albuquerque]}) Let $ E= E_1 \oplus E_2 $, where $ E $ is a real Banach space and $ E_1 $ is finite dimensional. Suppose that $ J \in C^{1}(E, \R) $ is even, $ J(0)=0 $, and that it verifies
\begin{itemize}
\item[$ (J_1) $] there exists $ \tau , r > 0 $ such that $ J(u) \geq \tau $ if $ ||u||_E = r $, $ u\in E_2 $,

\item[$ (J_2) $] there exists a finite-dimensional subspace $ \F \subset E $, with $ \dim E_1 < \dim \F $, and a constant $ \B > 0 $ such that $ \max\limits_{u\in \F} J(u) \leq \B $,

\item[$ (J_3) $] $ J $ satisfies the $ (PS)_c $ condition for all $ c\in (0, \B) $. 
\end{itemize}
Then, $ J $ possess at least $ \dim \F - \dim E_1 $ pairs of nontrivial critical points.  
\et

One should observe that, under conditions $ (f_1 ')-(f_4) $ and $ (g_1)-(g_5) $, we already have that $ I \in C^{1}(X, \R) $, is even, $ I(0)=0 $ and, from Lemma \ref{l11}, $ I $ verifies $ (J_1) $. It remains to prove that $ I $ also verifies $ (J_2) $ and $ (J_3) $.

Let $ k \in \N $, arbitrary but fixed, and $ Z \subset X $ a subspace of dimension $ k $, with norm $ ||\cdot||_Z $, which can be constructed by standard arguments.

\bl\label{l17}
There exists $ R>0 $ such that $ I(u) \leq 0 $ for all $ ||u||_Z \geq R $.
\el 
\begin{proof}
Since $ \dim Z < + \infty $ and from condition $ (f_4) $ and Lemma \ref{l2}-(b), we have
$$
I(u) \leq C_1 ||u||_{Z}^{p}+C_2 ||u||_{Z}^{N}+ C_3 ||u||_{Z}^{2p} + C_4 ||u||_{Z}^{p+ N} + C_5 ||u||_{Z}^{p+ \sigma} + C_6 ||u||_{Z}^{\sigma + N} - C_7 ||u||_{Z}^{q} \RA - \infty ,
$$
as $ ||u||_Z \RA + \infty $, since $ q > 2N $.
\end{proof}

\bl\label{l18}
There exists $ \eta > 0 $, sufficiently small, such that $ \max\limits_{u\in Z}I(u) \leq \eta $ and $ r_1 \alpha \left(\frac{1}{N}\right)^{\frac{1}{N-1}}\eta^{\frac{1}{N-1}}< \alpha_N $. 
\el
\begin{proof}
Let $ u\in Z \setminus \{0\} $. Then, from $ \dim Z < + \infty $, condition $ (f_4) $, Lemma \ref{l2}-(b) and Lemma \ref{l17}, one can find constants $ a_1=a_1(p, N, \sigma), a_2 = a_2(p, N, \sigma) >0 $ such that
$$
I(u) \leq a_1 ||u||_{Z}^{p}+ a_2 ||u||_{Z}^{\sigma + N} - 2 C_q M_q ||u||_{Z}^{q},
$$
where $||u||_q \geq M_q ||u||_Z $. Thus, we obtain a constant  $ a_3=a_3(p, N, \sigma)>0$
$$
\max\limits_{u\in Z} I(u) \leq \dfrac{a_3}{(C_q M_q)^{\frac{p}{q-p}}} .
$$
Therefore, taking $ C_q > 0 $ sufficiently large we find a value $ \eta >0 $ sufficiently small satisfying the desired conditions. 
\end{proof}

In the next lemma we guarantee that $ I $ satisfies the $ (PSC)_d $ condition for all $ d \in (0, \eta) $. Once the proof is very similar to the proof of Proposition \ref{p2}, we omit it here. We highlight that the validity of following lemma is possible only in virtue of Lemma \ref{lkey}.

\bl\label{l19}
The functional $ I $ satisfies condition $ (PSC)_d $ for all $ d\in (0, \eta) $.
\el

\begin{proof}[Proof of Theorem \ref{t2}]
From Lemmas \ref{l11}, \ref{l18} and \ref{l19} and an immediate application of Theorem \ref{t3}, with $ E=X $, $ E_1 = \{0\} $, $ \F = Z $, $ J=I $, $ \tau = m_\rho $, $ r= \rho $ and $ \B = \eta $, we get that $ I $ possess at least $ k $ nontrivial critical points. Therefore, as we can make $ k $ as large as we want, we conclude that (\ref{P}) has infinitely many solutions.
\end{proof}

\section{Related Problems}

In this final section we briefly discuss some similar problems that can be solved by these approach. One can consider the following equations 
\begin{equation} \label{P_1}
-\Delta_p u -\Delta_N u + |u|^{p-2}u + |u|^{N-2}u + \lambda (\ln|\cdot|\ast |u|^N)|u|^{N-2}u = f(u) \textrm{ \ in \ } \mathbb{R}^N,
\end{equation}
for $ 2 \leq p < N $ and 
\begin{equation} \label{P_2}
-\Delta_p u -\Delta_N u + a|u|^{p-2}u + b|u|^{N-2}u + \lambda (\ln|\cdot|\ast |u|^p)|u|^{p-2}u = f(u) \textrm{ \ in \ } \mathbb{R}^N,
\end{equation}
for $ \max\{2, \frac{N}{2}\}< p < N $. The reader should observe that the boundedness for $ p $ in this case is needed in order to obtain the equivalent of Lemma \ref{l13}. It is interesting that, in the particular case of  equations (\ref{P_1}) and (\ref{P_2}), one can obtain multiplicity results using genus theory, as in \cite{[boer3]}, which does not work for the general case (\ref{P}) since, even though considering additional growth conditions over $ g $, we do not have the desired geometry.    

On the other hand, in light of the same mentioned lemma, Lemma \ref{l13}, one can observe that this approach cannot be applied to a equation of the form
\begin{equation} \label{P_3}
-\Delta_p u -\Delta_N u + |u|^{p-2}u + |u|^{N-2}u + \lambda (\ln|\cdot|\ast |u|^N)|u|^{N-2}u + \gamma (\ln|\cdot|\ast |u|^p)|u|^{p-2}u = f(u) \textrm{ \ in \ } \mathbb{R}^N,
\end{equation}
since, in this case, we are not able to get rid of the logarithm term. 

Finally, we can solve versions of equations (\ref{P}), (\ref{P_1}) and (\ref{P_2}) considering two continuous potentials $ a, b : \RN \RA \R $ satisfying the following conditions

\noindent $ (c_0) \ a, b: \RN \RA \R \mbox{ are continuous, } \mathbb{Z}^N\mbox{-periodic , } a, b \in L^{\infty}(\RN), $
$$ \inf\limits_{x\in \RN} a(x) = a_0 > 0 \mbox{ and } \inf\limits_{x\in \RN} b(x) = b_0 > 0 .$$
One also one could try investigate the case where one, or both, the potentials $ a, b $ are not invariant under $ \mathbb{Z}^N $ translations but are asymptotically $ \mathbb{Z}^N $-periodic functions, that is, there exists a $ \mathbb{Z}^N $-periodic potential $ a_p : \RN \RA \R $ such that $ a_p $ satisfies $ (c_0) $,
$$
0 < \inf\limits_{x\in \RN} a(x) \leq a(x) \leq a_p(x) \ , \ \forall \ x \in \RN \leqno{(c_1)}
$$
and
$$
\lim\limits_{|x|\RA + \infty}|a(x) - a_p(x)| = 0 . \leqno{(c_2)}.
$$
The same for $ b $. To do that, one can combine the ideas presented here and in \cite{[alves2]}.


\begin{thebibliography}{2}

\bibitem{[albuquerque]} Albuquerque, F. S. B. (2014) Nonlinear Schrodinger elliptic systems involving exponential critical growth in $\mathbb{R}^2$, \textit{Electronic Journal of Differential Equations}. Vol. 2014, n. 59, pp. 1-12.

\bibitem{[pq4]} Alves, C. O., Ambrosio, V. and Isernia, T. (2019) Existence, multiplicity and concentration for a class of fractional $ p \&  q $ Laplacian problems in $ \mathbb{R}^{N} $, \textit{Communications on Pure \& Applied Analysis}. 18, 2009-2045.

\bibitem{[alves2]} Alves, C.O., Carrião, P.C.  and Miyagaki, O.H. (2001) Nonlinear Perturbations of a Periodic Elliptic Problem with Critical Growth, \textit{Journal of Mathematical Analysis and Applications}. 260, 133–146.

\bibitem{[alves]} Alves, C.O. and Figueiredo, G.M. (2019) Existence of positive solution for a planar Schrödinger-Poisson system with exponential growth, \textit{Journal of Mathematical Physics}. 60.

\bibitem{[ambrosseti]} Ambrosetti, A. and Rabinowitz P.H. (1973) Dual variational methods in critical point theory and applications, \textit{Journal of Functional Analysis}. 14, 349–381.

\bibitem{[pq3]} Barile, S. and Figueiredo, G. M. (2015) Existence of least energy positive, negative and nodal solutions for a class of $p \& q$-problems with potentials vanishing at infinity, \textit{Journal of Mathematical Analysis and Applications}. 427, 1205–1233.

\bibitem{[bartolo]} Bartolo, P., Benci, V., Fortunato, D. (1983) Abstract critical point theorems and applications to some nonlinear problems with “strong” resonance at infinity, \textit{Nonlinear Analysis: Theory, Methods \& Applications}. 7, 981–1012.



\bibitem{[boer3]} Böer, E. de S. and Miyagaki, O. H. (2021) Existence and multiplicity of solutions for the fractional $p$-Laplacian Choquard logarithmic equation involving a nonlinearity with exponential critical and subcritical growth, \textit{J. Math. Phys.} 62, 051507.

\bibitem{[4]} Bonheure, D., Cingolani, S., and Van Schaftingen, J. (2017) The logarithmic Choquard
equation: Sharp asymptotics and nondegeneracy of the groundstate. \textit{Journal of Functional
Analysis}, 272, 5255–5281.


\bibitem{[cjj]} Cingolani, S. and Jeanjean, L.  (2019) Stationary waves with prescribed $L^2$-norm for the planar Schr\"odinger-Poisson system, \textit{SIAM J. Math. Anal.} 51, 3533–3568.

\bibitem{[6]} Cingolani, S. and Weth, T. (2016) On the planar Schrödinger–Poisson system. \textit{Annales de l’Institut Henri Poincare (C) Non Linear Analysis}, 33, 169–197.




\bibitem{[doO]} do Ó, J. M. $N$-Laplacian equations in $\mathbb{R}^N$ with critical growth, \textit{Abstract and Applied Analysis}. 2 (1997) 301–315.


\bibitem{[10]} Du, M. and Weth, T. (2017) Ground states and high energy solutions of the planar Schrödinger–Poisson system. \textit{Nonlinearity}, 30, 3492–3515.

\bibitem{[pucci]} Fiscella, A.  and Pucci, P. (2019) $(p,N)$ equations with critical exponential nonlinearities in $\RN$, \textit{Journal of Mathematical Analysis and Applications}, 123379.

\bibitem{[12]} Fröhlich, H. (1954) Electrons in lattice fields. \textit{Advances in Physics}, 3, 325–361.
35

\bibitem{[pq1]} He, C. and Li, G. (2008) The regularity of weak solutions to nonlinear scalar field elliptic equations containing $p\& q$-Laplacians. \textit{Annales Academiae Scientiarum Fennicae}, 33, 337-371. 





\bibitem{[lam]} Lam, N. and Lu, G. (2012) Existence and multiplicity of solutions to equations of $N$-Laplacian type with critical exponential growth in $\RN$, \textit{Journal of Functional Analysis}. 262, 1132–1165.

\bibitem{[w6]} Li, Q. and Yang, Z. (2016) Multiple solutions for a class of fractional quasi-linear equations with critical exponential growth in $\RN $, \textit{Complex Variables and Elliptic Equations}. 61, 969–983.

\bibitem{[15]} Lieb, E. H. (1983) Sharp constants in the Hardy-Littlewood-Sobolev and related inequalities. \textit{The Annals of Mathematics}, 118, 349.

\bibitem{[lions]} Lions, P. L. (1984) The concentration-compactness principle in the Calculus of Variations. The Locally compact case, part 2, \textit{Annales de l’Institut Henri Poincaré C}, Analyse Non Linéaire. 1, 223–283.

\bibitem{[16]} Lions, P. L. (1987) Solutions of Hartree-Fock equations for Coulomb systems. \textit{Communications in Mathematical Physics}, 109, 33–97.

\bibitem{[17]} Moser, J. (1971) A Sharp form of an inequality by N. Trudinger. \textit{Indiana University Mathematics Journal}, 20, 1077–1092.

\bibitem{[18]} Penrose, R. (1996) On gravity’s role in quantum state reduction. \textit{General Relativity and Gravitation}, 28, 581–600.

\bibitem{[20]} Ruiz, D. (2006) The Schrödinger–Poisson equation under the effect of a nonlinear local
term. \textit{Journal of Functional Analysis}, 237, 655–674.

\bibitem{[silva]} Silva, E. A. B. (1988) \textit{Critical point theorems and applications to differential equations}, PhD. Thesis, University of Wisconsin-Madison.


\bibitem{[21]} Stubbe, J. (2008) Bound states of two-dimensional Schr\"{o}dinger-Newton equations.
arXiv:0807.4059 \textit{[math-ph]}, arXiv: 0807.4059.


\bibitem{[pq2]} Tanaka, M. (2014) Generalized eigenvalue problems for $( p , q )$ -Laplacian with indefinite weight, \textit{Journal of Mathematical Analysis and Applications}. 419, 1181–1192.

\bibitem{[22]} Wilson, A. J. C. (1955) Untersuchungen über die Elektronentheorie der Kristalle by S. I.
Pekar. \textit{Acta Crystallographica}, 8, 70–70.

\bibitem{[zhang]} Zhang, C. and Chen, L. Concentration-Compactness Principle of Singular Trudinger-Moser Inequalities in $\RN$ and $n$-Laplace Equations, \textit{Advanced Nonlinear Studies}. 18 (2018) 567–585.


\end{thebibliography}
\end{document}